\documentclass[12pt,fleqn]{article}
\usepackage{graphicx}

\usepackage{latexsym}

\usepackage{amsmath}
\usepackage{amsthm}
\usepackage{amssymb}
\usepackage{amsfonts}

\numberwithin{equation}{section}

\begin{document}

\newcommand{\rf}[1]{(\ref{#1})}
\newcommand{\rff}[2]{(\ref{#1}\ref{#2})}

\newcommand{\ba}{\begin{array}}
\newcommand{\ea}{\end{array}}

\newcommand{\be}{\begin{equation}}
\newcommand{\ee}{\end{equation}}

\newcommand{\const}{{\rm const}}
\newcommand{\ep}{\varepsilon}
\newcommand{\Cl}{{\cal C}}
\newcommand{\rr}{{\vec r}}
\newcommand{\ph}{\varphi}
\newcommand{\R}{{\mathbb R}}  
\newcommand{\N}{{\mathbb N}}
\newcommand{\Z}{{\mathbb Z}}

\newcommand{\e}{{\bf e}}

\newcommand{\m}{\left( \ba{r}}
\newcommand{\ema}{\ea \right)}
\newcommand{\mm}{\left( \ba{cc}}
\newcommand{\miv}{\left( \ba{cccc}}

\newcommand{\scal}[2]{\mbox{$\langle #1 \! \mid #2 \rangle $}}
\newcommand{\ods}{\par \vspace{0.3cm} \par}
\newcommand{\dis}{\displaystyle }
\newcommand{\mc}{\multicolumn}
\newcommand{\no}{\par \noindent}

\newcommand{\sinc}{ {\rm sinc\,} }
\newcommand{\tanhc}{{\rm tanhc} }
\newcommand{\tanc}{{\rm tanc} }

\newtheorem{prop}{Proposition}[section]
\newtheorem{Th}[prop]{Theorem}
\newtheorem{lem}[prop]{Lemma}
\newtheorem{rem}[prop]{Remark}
\newtheorem{cor}[prop]{Corollary}
\newtheorem{Def}[prop]{Definition}
\newtheorem{open}{Open problem}
\newtheorem{ex}[prop]{Example}
\newtheorem{exer}[prop]{Exercise}

\newenvironment{Proof}{\par \vspace{2ex} \par
\noindent \small {\it Proof:}}{\hfill $\Box$ 
\vspace{2ex} \par }

\title{\bf 
Locally exact modifications of  numerical schemes}
\author{
 {\bf Jan L.\ Cie\'sli\'nski}\thanks{\footnotesize
 e-mail: \tt janek\,@\,alpha.uwb.edu.pl}
\\ {\footnotesize Uniwersytet w Bia{\l}ymstoku,
Wydzia{\l} Fizyki}
\\ {\footnotesize ul.\ Lipowa 41, 15-424
Bia{\l}ystok, Poland}
}

\date{}

\maketitle

\begin{abstract}
We present a new class of exponential integrators for ordinary differential equations: locally exact modifications of   known  numerical schemes. Local exactness means that they preserve the linearization of the original system at every point.  In particular, locally exact integrators preserve all fixed points and are A-stable. We apply this approach to popular schemes including  Euler schemes, implicit midpoint rule and trapezoidal rule. We found locally exact modifications of discrete gradient schemes (for symmetric discrete gradients and coordinate increment discrete gradients) preserving their main geometric property: exact conservation of the energy integral (for arbitrary multidimensional Hamiltonian systems in canonical coordinates).  Numerical experiments for a 2-dimensional anharmonic oscillator show that locally exact schemes have very good accuracy in the neighbourhood of stable equilibrium, much higher than suggested by the order of new schemes (locally exact modification sometimes increases the order but in many cases leaves it unchanged). 
\end{abstract}

\ods




\section{Introduction}

The motivation for introducing ``locally exact'' discretizations is quite natural. Considering a numerical scheme for a dynamical system with periodic solutions (for instance: a nonlinear pendulum) we may ask whether the numerical scheme recovers the period of small oscillations. For a fixed finite time step the answer is usualy negative. However, many numerical scheme (perhaps all of them) admit modifications which preserve the period of small oscillations for a fixed (not necessarily small) time step $h$.  An unusual feature of our approach is that instead of taking the limit $h \rightarrow 0$, we consider the limit $x_n \approx \bar x$ (where $\bar x$ is the stable equilibrium).  
As the next step we consider a linearization around any fixed $\bar x$ (then, in order to preserve the condition $x_n \approx \bar x$, some evolution of $\bar x$ is necessary). 
In other words, we combine two known procedures: the approximation of nonlinear systems by linear equations and explicit exact discretizations of linear equations with constant coefficients. An essential  novelty consists in applying these procedures to a modified numerical scheme containing free functional parameters.  
  The case of small oscillations was presented in \cite{CR-long}. Our method works perfectly for discrete gradient schemes \cite{LaG,MQR2,QT}.  We succeded in modifying  discrete gradient schemes in a locally exact way without spoiling their main geometric property: the exact conservation of the energy integral. The  preservation of geometric properties by numerical algorithms is of considerable advantage \cite{HLW,Is}. Promising results on  one-dimensional Hamiltonian systems can be found in \cite{CR-PRE,CR-BIT} (by one-dimensional Hamiltonian system we mean a Hamiltonian system with one degree of freedom).  In this paper we extend our approach on the case of multidimensional canonical Hamiltonian systems. Moreover, we present locally exact modifications of forward and backward Euler schemes, implicit midpoint rule and trapezoidal rule.

A notion closely related to our local exactness has been proposed a long time ago \cite{Pope}, see also \cite{Law}.  Recently, similar concept appeared under the name of linearization-preserving preprocessing \cite{MQTse}, see also  below (section~\ref{sec-locex}). Our approach has  also some similarities with the Mickens approach \cite{Mic}, Gautschi-type methods \cite{HoL2} and, most of all, with the exponential integrators technique \cite{CCO,HoL1,MW}. The definition of an exponential integrator is so wide (e.g., ``a numerical method which involves an exponential function of the Jacobian'', \cite{HoL1}) that our  schemes can be considered as special exponential integrators. 
  In spite of some similarities and overlappings, our approach differs from all methods mentioned above. In particular, according to our best knowledge, discrete gradient schemes have never been treated or modified in the framework of exponential integrators and/or linearization-preserving preprocessing.

\section{Locally exact numerical schemes}

We consider an ordinary differential equation (ODE) with the general solution ${\pmb x} (t)$ (satisfying the initial condition ${\pmb x} (t_0) = {\pmb x}_0$), and a difference equation with the  general solution ${\pmb x}_n$. 
The difference equation 
 is the exact discretization of  the considered ODE  if \ ${\pmb x}_n = {\pmb x} (t_n)$. 

\subsection{Exact discretization of linear systems}
\label{sec-exact}

It is well known that any linear ODE with constant coefficients admits the exact discretization in an explicit form \cite{Potts}, see also \cite{Ag,CR-ade,Mic}. We summarize these results as follows, compare \cite{CR-BIT} (Theorem 3.1).

\begin{prop} \label{prop-exact} 
Any linear equation with constant coefficients, represented in the matrix form by
\be  \label{genlin}
  \frac{d  {\pmb x}}{d t} = A {\pmb x} + {\pmb b} \ ,
\ee
(where ${\pmb x} = {\pmb x} (t) \in \R^d$, ${\pmb b} = \const \in \R^d$ and $A$ is a constant invertible real $d\times d$ matrix) admits the exact discretization given by
\be   \label{genex}
   {\pmb x}_{n+1} - {\pmb x}_n  = (e^{h_n A} - I) A^{-1} \left( A {\pmb x}_n + {\pmb b} \right)   \ , 
\ee
where \ $h_n = t_{n+1} - t_n$ \  is the time step and $I$ is the identity matrix. 
\end{prop}
\ods

\begin{cor}  We may look at the exact discretization \rf{Del-ex} as a modification of the forward Euler scheme. Indeed,  we can rewrite 
 \rf{genex} as 
\be  \label{Del-ex}
  \Delta_n^{-1} \left(  {\pmb x}_{n+1} - {\pmb x}_n  \right) = A {\pmb x}_n + {\pmb b} \ ,
\ee
where  $\Delta_n$, defined by  
\be  
   \Delta_n =  A^{-1} (e^{h_n A} - I) \ , 
\ee
is  a matrix ``perturbation'' of the time step $h_n$. Indeed,  $ \Delta_n = h_n I + O (h_n^2)$. 
\end{cor}
 
Exact discretizations found direct application in numerical  treatment of the classical Kepler problem \cite{Ci-Kep,Ci-Koz,Ci-oscyl}. The exact discretization of the harmonic oscillator equation can also be used to in the integration of some partial differential equations by Fourier transformation \cite{Ci-oscyl}.  

The central topic of our paper is another fruitful direction of using exact integrators,  namely the so called locally exact discretizations \cite{Ci-oscyl,CR-PRE}.

\subsection{Local exactness}
\label{sec-locex}

Motivated by the results of \cite{CR-long,CR-PRE} we propose the following definition.

\begin{Def} \label{def-locex-at}
A numerical scheme ${\pmb x}_{n+1} = \Phi ({\pmb x}_n, h_n)$ for an autonomous equation $\dot {\pmb x} = F ({\pmb x})$  is {\it locally exact} at $\pmb{\bar x}$ if its linearization around $\pmb{\bar x}$ is identical with the exact discretization of the differential equation linearized around $\pmb{\bar x}$. 
\end{Def}

 The simplest choice is $\bar x = x_0$, where $V'(x_0)=0$  (small oscillations around the stable equilibrium). In this case $\delta_n$ does not depend on $n$. The resulting scheme, known as MOD-GR (compare \cite{CR-BIT,CR-CPC}),  was first presented in \cite{CR-long}. In  \cite{CR-PRE} we considered the case 
 $\bar x = x_n$ (GR-LEX) and its symmetric (time-reversible) modification \ $\bar x = \frac{1}{2} (x_n + x_{n+1})$ (GR-SLEX). Note that in both cases  $\bar x$  changes at every step. 
The scheme MOD-GR is locally exact at a stable equilibrium only, GR-LEX is locally exact at $x_n$ (for any $n$), and, finally, GR-SLEX is locally exact at $\frac{1}{2} (x_n + x_{n+1})$ (for any $n$). In the case of implicit numerical schemes the function $\Phi ({\pmb x}_n, h_n)$ is, of course, an implicit function.

\begin{Def}  \label{def-locex}
A numerical scheme ${\pmb x}_{n+1} = \Phi ({\pmb x}_n, h_n)$ for an autonomous equation $\dot {\pmb x} = F ({\pmb x})$  is {\it locally exact} if there exists a sequence $\pmb{\bar x}_n$ such that $\pmb{\bar x}_n  - {\pmb x}_n = O (h_n)$ and the scheme is locally exact at $\pmb{\bar x}_n $ (for any $n$). 
\end{Def}

Therefore GR-LEX and GR-SLEX are locally exact. 
Similar concept (``linearization-preserving'' schemes)  has been recently formulated in \cite{MQTse} (Definition 2.1). An integrator is said to be linearization-preserving if it is linearization preserving at all fixed points. 
This definition is weaker than our Definition~\ref{def-locex}. All locally exact schemes are linearization-preserving but, in general, linearization-preserving schemes do not have to be   locally exact in our sense. 
For instance, the scheme MOD-GR (see \cite{CR-long,CR-BIT}) is linearization-preserving (provided that $V (x)$ has only one stable equilibrium) but is not locally exact.  The problem of finding the best sequence $\pmb{\bar x}_n$  for a given numerical scheme seems to be interesting but has not been considered yet.

\subsection{Exact discretization of linearized equations}
\label{sec-ex-lin}

As an immediate consequence of Definition~\ref{def-locex} we have that the exact discretization of the linearization of a given nonlinear system is locally exact. 
We confine ourselves to autonomous systems of the form
\be  \label{firstord}
  \pmb{\dot x} = F ({\pmb x}) 
\ee
where ${\pmb x} (t) \in \R^d$. 
If ${\pmb x}$ is near a fixed $\pmb{\bar x}$, then \rf{firstord} can be approximated by
\be  \label{approx1}
\dot {\pmb \xi} = F' (\pmb{\bar x}) {\pmb \xi}  +  F (\pmb{\bar x}) 
\ee
where $F'$ is the Fr\'echet derivative (Jacobi matrix) of $F$ and 
\be  \label{xi}
 {\pmb \xi} = {\pmb x} - \pmb{\bar x} \ . 
\ee
The exact discretization of the approximated equation \rf{approx1} is given by
\be  \label{exact-approx}
{\pmb \xi}_{n+1} = e^{ h_n F' (\pmb{\bar x})} {\pmb \xi}_n + \left( e^{ h_n F' (\pmb{\bar x})} - I \right) \left( F' (\pmb{\bar x}) \right)^{-1} F (\pmb{\bar x}) \ , 
\ee
provided that \ $\det F' (\pmb{\bar x}) \neq 0$. This assumption, however, is not necessary because function
\be
\ph_1 (x) := \frac{e^x - 1}{x}    \qquad  (\text{and} \ \ph_1 (0):= 1) 
\ee
is analytic also at $x=0$. 
\ods
\begin{prop}  \label{cor-linex}
The exact discretization of the linearization of $\pmb{\dot x} = F (\pmb{x})$ around $\pmb{\bar x}$ is given by: 
\be  \label{lin-ex}
{\pmb x}_{n+1} - \pmb{\bar x} = e^{ h_n F' (\pmb{\bar x})} \left( x_n - \pmb{\bar x} \right) + h_n \ph_1 (h_n F' (\pmb{\bar x}))  F (\pmb{\bar x})  \ . 
\ee
The scheme \rf{lin-ex} is locally exact. 
\end{prop}
\ods

If we put $\pmb{\bar x} = {\pmb x}_n$ into \rf{lin-ex} (thus changing $\pmb{\bar x}$ at each step), then ${\pmb \xi}_n = 0$ and we obtain:
\be  \label{le-pred}
  {\pmb x}_{n+1} = {\pmb x}_n + \left( e^{  h_n F' ({\pmb x}_n) } - I \right) ( F'({\pmb x}_n) )^{-1} F ({\pmb x}_n) 
\ee
This method is well known as the exponential difference equation (see \cite{Pope}, formula (3.4)) or, more recently, as the exponential Euler method \cite{MW}.

\section{Locally exact modifications of popular one-step numerical schemes }
\label{sec-pop}

In order to illustrate the concept of locally exact modifications, we apply this procedure to  the following class of numerical methods:
\be \label{scheme-class}
  {\pmb y}_{n+1} - {\pmb y}_n =h   \Psi ({\pmb y}_n, {\pmb y}_{n+1} ) \  , 
\ee
where $\Psi (\pmb{y}, \pmb{y}) = F (\pmb{y})$  and  \ $\dot {\pmb y} = F ({\pmb y})$ \ is equation to be solved. This class   contains, among others, the following numerical schemes: \par
\renewcommand{\arraystretch}{1.5}\par
\begin{tabular}{lll}
$\bullet$ &  \text{explicit Euler scheme (EEU):} & $\Psi ({\pmb y}_n, {\pmb y}_{n+1} ) = F (\pmb{y}_n)$  , \\ 
$\bullet$ & \text{implicit Euler scheme (IEU):} & $\Psi ({\pmb y}_n, {\pmb y}_{n+1} ) = F (\pmb{y}_{n+1})$  , \\ 
$\bullet$ & \text{implicit midpoint rule (IMP):} & $\Psi ({\pmb y}_n, {\pmb y}_{n+1} ) = F (\frac{ \pmb{y}_{n} + \pmb{y}_{n+1}}{2} )$  , \\
$\bullet$ & \text{trapezoidal rule (TR):} & $\Psi ({\pmb y}_n, {\pmb y}_{n+1} ) = \frac{1}{2} \left( F (\pmb{y}_n) + F(\pmb{y}_{n+1}) \right)$ .
\end{tabular} \par
\renewcommand{\arraystretch}{1} \ods

A natural  non-standard modification of \rf{scheme-class} can be obtained by replacing $h$ by a matrix  $\pmb{\delta}$. Thus we consider the following familiy of modified numerical schemes:   
\be \label{scheme-gen}
  {\pmb y}_{n+1} - {\pmb y}_n = \pmb{\delta} ({\pmb{\bar y}_n},h_n) \Psi ({\pmb y}_n, {\pmb y}_{n+1} ) \ ,
\ee
where $\pmb{\bar y}_n$ is described in Definition~\ref{def-locex} and $\pmb{\delta} ({\pmb{\bar y}_n},h_n)$ is a matrix-valued function. 
We assume  consistency conditions:
\be  \label{consis}
\lim_{h \rightarrow 0} \frac{ \pmb{\delta} ({\pmb{\bar y}_n},h) }{h} = I  , \quad \Psi ({\pmb y}, {\pmb y}) = F ({\pmb y})  , \quad  \Psi_1 ({\pmb y}, {\pmb y}) + \Psi_2 ({\pmb y}, {\pmb y})    = F'({\pmb y})  , 
\ee
where $I$ is $m\times m$ identity matrix, $\Psi_1, \Psi_2$ are partial Fr\'echet derivative of $\Psi$ with respect to the first and second vector variable, respectively (thus $\Psi_1, \Psi_2$ are matrices). We also denote
\be  \label{xx}
 {\bar \Psi} = \Psi ({\pmb{\bar y}_n}, {\pmb{\bar y}_n}) \ , \quad {\bar \Psi}_1 =  \Psi_1 ({ \pmb{\bar y}_n}, {\pmb{\bar y}_n}) \ , \quad {\bar \Psi}_2 =  \Psi_2 ({ \pmb{\bar y}_n}, {\pmb{\bar y}_n}) \ .
\ee

\begin{Th}  \label{prop-legen}
Numerical scheme \rf{scheme-gen}, where $\pmb{\bar y}_n$ is any sequence described in Definition~\ref{def-locex} and $\Psi$ satisfies \rf{consis}, is locally exact for
\be    \label{delta-gen}
\pmb{\delta} ({\pmb{\bar y}_n},h_n) = h_n  \varphi_1 (h_n F' (\pmb{\bar y}_n))   \left( I + h_n {\bar \Psi}_2  \varphi_1 (h_n F' (\pmb{\bar y}_n))   \right)^{-1} \  .
\ee
\end{Th}

\begin{Proof} The exact discretization of the linearization of equation ${\pmb{\dot y}} = F ({\pmb y})$ is given by  \rf{exact-approx}. The linearization of the scheme \rf{scheme-gen} (at $\pmb{\bar y}_n$)  reads
\be  \label{ksi1}
  {\pmb \nu}_{n+1} - {\pmb  \nu}_n = \pmb{\delta} ( {\bar \Psi}_1 {\pmb \nu}_n + {\bar \Psi}_2 {\pmb \nu}_{n+1} +  {\bar \Psi}) \ , 
\ee
where ${\pmb \nu}_n =  \pmb{y}_n - \pmb{\bar y}_n$, \  $\pmb{\delta} = \pmb{\delta} (\pmb{\bar y}_n, h_n)$ and we use \rf{xx}. Identifying \rf{ksi1} with \rf{exact-approx} and assuming invertibility of  \ $I - \pmb{\delta} {\bar \Psi}_2$ \  we get a system of two equations: 
\be \ba{l}  \label{fir}
(I - \pmb{\delta} {\bar \Psi}_2)^{-1} (I + \pmb{\delta} {\bar \Psi}_1) = e^{h_n F'} \ ,   \\[1em]
(I - \pmb{\delta} {\bar \Psi}_2)^{-1} \pmb{\delta} {\bar \Psi} = h_n  \varphi_1 (h_n F') F \ ,
\ea \ee
where $F = F (\pmb {\bar y}_n)$ and  $F'= F' (\pmb {\bar y}_n)$.  
Equations \rf{fir} imply 
\be \ba{l}  \label{symm}
\pmb{\delta} \left( {\bar \Psi}_1 + {\bar \Psi}_2  e^{h_n F'} \right) = e^{h_n F'} - I  \ , \\[1em]
\pmb{\delta} \left(   {\bar \Psi}  + h_n  {\bar \Psi}_2   \varphi_1 (h_n F') F  \right) = h_n  \varphi_1 (h_n F') F \ .
\ea \ee
Using \rf{consis} and \rf{xx} (i.e., $\bar \Psi = F$, \ $\bar{\Psi}_1 + \bar{\Psi}_2 = F'$), and replacing $ e^{h_n F'}$  by $I + h_n F'  \varphi_1 (h_n F')$, we rewrite \rf{symm} as follows
\be  \ba{l}  \label{both}
\pmb{\delta} \left(  I + h_n \bar{\Psi}_2  \varphi_1 (h_n F') \right) F' = h_n \varphi_1 (h_n F') F' \ , \\[1em]
\pmb{\delta} \left(  I + h_n \bar{\Psi}_2  \varphi_1 (h_n F') \right) F = h_n \varphi_1 (h_n F') F \ . 
\ea \ee
One can easily see that $\pmb{\delta}$ given by  \rf{delta-gen} satisfies simultaneously both equations \rf{both} (actually \rf{delta-gen} is also necessary provided that $F'$ is invertible).
\end{Proof}

\ods

By straightforward calculation we  can express  matrix  $\pmb{\delta}$ in the following (equivalent) form: 
\be  \label{del-tanhc}
   \pmb{\delta} ({\pmb{\bar y}_n}, h_n) =  h_n  \tanhc \frac{h_n F'}{2}  \left( I + \frac{1}{2} h_n \left( {\bar \Psi}_2 - {\bar \Psi}_1 \right) \tanhc \frac{h_n F'}{2}  \right)^{-1}   ,
\ee
where $\dis   \tanhc (z) \equiv  z^{-1} \tanh (z)$ \ 
is analytic at $z = 0$. For small $h_n$ we have
\[
      \pmb{\delta} ({\pmb{\bar y}_n}, h_n)  = h_n I + \frac{1}{2}  h_n^2 (\bar \Psi_1 - \bar \Psi_2) + O (h_n^3) \ .
\]
Therefore $\pmb{\delta}$ exists for sufficiently small $h_n$. In other words, numerical scheme \rf{scheme-gen} admits a locally exact modification for sufficiently small $h_n$.

\subsection*{Locally exact explicit Euler scheme}

Specializing $\Psi ({\pmb y}_n, {\pmb y}_{n+1} ) = F (\pmb{y}_n)$ 
we obtain a class of locally exact modifications of the explicit Euler scheme 
\be  \label{lex-explEuler}
{\pmb x}_{n+1} = {\pmb x}_n + \big( e^{ h_n F' (\pmb{\bar x})} - I \big) ( F' (\pmb{\bar x}) )^{-1} F ({\pmb x}_n) \ .
\ee
In this case the choice $\pmb{\bar x} = {\pmb x}_n$ seems to be most natural.

\begin{prop}
Locally exact modification of the explicit Euler scheme (EEU-LEX) is given by
\be  \label{lee-Euler}
{\pmb x}_{n+1} = {\pmb x}_n + \big( e^{ h_n F' ({\pmb x}_n)} - I \big) ( F' ({\pmb x}_n) )^{-1} F ({\pmb x}_n) \ .
\ee
\end{prop}

\no Note that   \rf{lee-Euler} coincides with \rf{le-pred}).

\subsection*{Locally exact implicit Euler schemes}

For $\Psi ({\pmb y}_n, {\pmb y}_{n+1} ) = F (\pmb{y}_{n+1})$ we get
\be  \label{lex-impEuler-gen}
{\pmb x}_{n+1} = {\pmb x}_n + \big(  I - e^{-  h_n F' (\pmb{\bar x})} \big) ( F' (\pmb{\bar x}) )^{-1} F ({\pmb x}_{n+1}) \ ,
\ee
or, in an equivalent way,
\be
{\pmb x}_{n+1} = {\pmb x}_n +  h_n \ \ph_1 ( -  h_n F'(\pmb{\bar x}) ) \ F ({\pmb x}_{n+1}) \ .
\ee
In this case it is not clear what is the most natural identification. We can choose either $\pmb{\bar x} = {\pmb x}_n$, or $\pmb{\bar x} = {\pmb x}_{n+1}$.  

\begin{prop}
Locally exact modification of the implicit Euler scheme is given either by
\be  \label{leie-Euler}
{\pmb x}_{n+1} = {\pmb x}_n + \big(  I - e^{-  h_n F' ({\pmb x}_n)} \big) ( F' ({\pmb x}_n) )^{-1} F ({\pmb x}_{n+1}) \ ,
\ee
which will be called IEU-LEX, or  by
\be
\label{leii-Euler}
{\pmb x}_{n+1} = {\pmb x}_n + \big(  I - e^{-  h_n F' ({\pmb x}_{n+1})} \big) ( F' ({\pmb x}_{n+1}) )^{-1} F ({\pmb x}_{n+1}) \  ,
\ee
called IEU-ILEX. 
\end{prop}
Both locally exact implicit Euler schemes are of second order. IEU-ILEX is a little bit more accurate, but IEU-LEX is more effective when computational costs are taken into account (see Section~\ref{sec-exp}). 

\subsection*{Locally exact implicit midpoint rule} 

We take $\Psi ({\pmb y}_n, {\pmb y}_{n+1} ) = F (\frac{ \pmb{y}_{n} + \pmb{y}_{n+1}}{2} )$. Then, using \rf{del-tanhc}, we get
\be  \label{imr-deltan}
\pmb{\delta}_n  = 2 ( F' (\pmb{\bar x}) )^{-1} \tanh \left( \frac{h_n F'(\pmb{\bar x})}{2} \right)   . 
\ee
Thus locally exact modification of the implicit midpoint rule is given by
\be  \label{le-midpoint}
{\pmb x}_{n+1} - {\pmb x}_n  = h_n \left( \tanhc \frac{ h_n F'(\pmb{\bar x})}{2} \right) \ F \left( \frac{{\pmb x}_{n+1} + {\pmb x}_n}{2} \right) 
\ee
The most natural choice of $\pmb{\bar x}$ seems to be at the midpoint, $\pmb{\bar x} = \frac{1}{2} \left( {\pmb x}_n + {\pmb x}_{n+1} \right)$. The obtained scheme will be called IMP-SLEX. However, in order to diminish the computation cost, the choice $\pmb{\bar x} = {\pmb x}_n$, i.e., IMP-LEX,  can also be considered (because then Jacobian $F' (\pmb{\bar x})$ is evaluated outside iteration loops).

\subsection*{Locally exact trapezoidal rule} 

Taking into account $\Psi ({\pmb y}_n, {\pmb y}_{n+1} ) = \frac{1}{2} \left( F (\pmb{y}_n) + F(\pmb{y}_{n+1}) \right)$ and \rf{del-tanhc}, we get \rf{imr-deltan}, as well. 
Thus locally exact modification of the trapezoidal rule reads
\be  \label{le-trapez}
{\pmb x}_{n+1} - {\pmb x}_n  =  h_n \left( \tanhc \frac{ h_n F'(\pmb{\bar x})}{2} \right) \ \frac{ F ({\pmb x}_{n+1}) + F ({\pmb x}_n) }{2} \ . 
\ee
There are two natural choices of $\pmb{\bar x}$. Either (in order  to minimize the computational costs) we can take $\pmb{\bar x} = {\pmb x}_n$, obtaining  TR-LEX,  or (in order to get   time-reversible scheme  TR-SLEX) we can take $\pmb{\bar x} = \frac{1}{2} \left( {\pmb x}_n + {\pmb x}_{n+1} \right)$.

\section{Locally exact discrete gradient schemes for multidimensional Hamiltonian systems}
\label{sec-locex-multi}

In this section we construct   
energy-preserving locally exact discrete gradient schemes for arbitrary multidimensional Hamiltonian systems in canonical coordinates:
\be \label{multiham}
  {\dot x}^k = \frac{\partial H}{\partial p^k} \ , \quad {\dot p}^k = - \frac{\partial H}{\partial x^k} \ .
\ee
where $k = 1,\ldots,m$. 
We obtain two different numerical schemes using either  symmetric discrete gradient or coordinate increment discrete gradient. 

\subsection{Linearization of Hamiltonian systems}

We denote $\pmb{x} = (x^1,\ldots, x^m)^T$, $\pmb{p} = (p^1,\ldots, p^m)$, etc. 
The linearization of \rf{multiham} around \ $\pmb{\bar x}, \pmb{\bar p}$ \ is given by:
\be \ba{l}  \label{lin-mham}
 \pmb{\dot \xi} = H_{\pmb{p}} + H_{\pmb{p} \pmb{x}} \pmb{\xi} + H_{\pmb{p} \pmb{p}} \pmb{\eta} \ , \\[2ex]
\pmb{\dot \eta} = - H_{\pmb{x}} - H_{\pmb{x} \pmb{x}} \pmb{\xi} - H_{\pmb{x} \pmb{p}} \pmb{\eta} \ , 
\ea \ee
where $\pmb{\xi} = \pmb{x}-\pmb{\bar x}$, \ $\pmb{\eta} = \pmb{p}-\pmb{\bar p}$, \  $H_{\pmb{p}}$ is a vector with components $\frac{\partial H}{\partial p^k}$,   $H_{\pmb{p} \pmb{x}}$ is a matrix with coefficients $\frac{\partial^2 H}{\partial p^j \partial x^k}$, etc., and derivatives of $H$ are evaluated at $\pmb{\bar x}, \pmb{\bar p}$.  
Note that $m\times m$ matrices $H_{\pmb{x} \pmb{x}}$, $H_{\pmb{x} \pmb{p}}$, $H_{\pmb{p} \pmb{x}}$, $H_{\pmb{p} \pmb{p}}$ satisfy 
\be
 H_{\pmb{p} \pmb{p}}^T = H_{\pmb{p} \pmb{p}} \ , \quad H_{\pmb{x} \pmb{x}}^T = H_{\pmb{x} \pmb{x}}  \ , \quad H_{\pmb{x} \pmb{p}}^T = H_{\pmb{p} \pmb{x}}  
\ee
Equations \rf{lin-mham} can be rewritten also as
\be \label{lin-mmh}
\frac{d}{d t} \m \pmb{\xi} \\ \pmb{\eta} \ema = F' \m \pmb{\xi} \\ \pmb{\eta} \ema  + F 
\ee
where 
\be  \label{F-F'}
 F = \m H_{\pmb{p}} \\ - H_{\pmb{x}} \ema \ , \quad F' = \mm H_{\pmb{p}\pmb{x}} & H_{\pmb{p}\pmb{p}} \\ - H_{\pmb{x}\pmb{x}} & - H_{\pmb{x}\pmb{p}} \ema \ .
\ee

\begin{cor}  \label{cor-linham} 
The exact discretization of linearized Hamiltonian equations \rf{lin-mmh} is given by 
\be  \label{lin-hex}
 \m \pmb{\xi}_{n+1} \\ \pmb{\eta}_{n+1} \ema = e^{h_n F'} \m \pmb{\xi}_n \\ \pmb{\eta}_n \ema + h_n \ph_1 (h_n F')  F \ .
\ee
\end{cor}

\begin{Proof}
 The exact discretization of \rf{lin-mmh} is given by \rf{lin-hex} which follows immediately from Proposition~\ref{cor-linex}. 
\end{Proof}

\subsection{Discrete gradients in the multidimensional case} 

Considering multidimensional Hamiltonian systems we will denote 
\be  \label{y=xp}
\pmb{y} := \m \pmb{x} \\  \pmb{p} \ema  , \quad   \pmb{y}_n := \m \pmb{x}_n \\ \pmb{p}_n \ema ,  
\ee
where $\pmb{x} \in \R^m$, $\pmb{p} \in \R^m$, $\pmb{y} \in \R^{2 m}$, etc. In other words, 
\[ 
y^1 = x^1  , \ \  y^2 = x^2 , \ \ldots \ y^m = x^m , \ \  y^{m+1} = p^1  ,\ \ldots \ y^{2m} = p^m  . 
\]
A discrete gradient, denoted by $\bar \nabla H (\pmb{y}_n, \pmb{y}_{n+1})$,  
\be  \label{barnablaH}
 {\bar \nabla} H   = \left( \frac{ {\Delta} H}{\Delta y^1}, \frac{ {\Delta} H}{\Delta y^2}, \ldots,\frac{ {\Delta} H}{\Delta y^{2 m}}  \right) \equiv \left( \frac{ \Delta H}{\Delta \pmb{x}}, \frac{ \Delta H }{\Delta \pmb{p}} \right) 
\ee
is defined as an $\R^{2 m}$-valued  function of $\pmb{y}_n,  \pmb{y}_{n+1}$ such that \cite{Gon,MQR2}: 
\be \ba{l}\dis \label{disgrad-def}
 \sum_{k=1}^{2 m}   \frac{ {\Delta} H }{\Delta y^k } \left( y_{n+1}^k - y_n^k  \right)  = H (\pmb{y}_{n+1}) - H (\pmb{y}_n) \ , 
\\[4ex]\dis
{\bar \nabla} H (\pmb{y}, \pmb{y})  = \left( H_{\pmb{x}}, H_{\pmb{p}} \right) \ , 
\ea \ee
where $H_{\pmb{x}}, H_{\pmb{p}}$ are evaluated at $\pmb{y} = (\pmb{x}, \pmb{p})$ and we define
\be  \label{grad-lim}
{\bar \nabla} H (\pmb{y}, \pmb{y}) = \lim_{\pmb{\tilde y} \rightarrow \pmb{y}} {\bar \nabla} H (\pmb{y}, \pmb{\tilde y})
\ee 
when necessary. It is well known (see, e.g., \cite{Gon,Gree,LaG,MQR2}) that any discrete gradient defines an energy-preserving numerical scheme 
\be  \label{disgrad}
   \pmb{y}_{n+1} = \pmb{y}_n + h  \bar \nabla H \ . 
\ee

Discrete gradients are non-unique, compare \cite{Gon,IA,MQR2}. Here we confine ourselves to the simplest discrete gradients. We consider coordinate increment discrete gradient, first proposed by Itoh and Abe \cite{IA}, and its symmetrization (see below).  The coordinate increment discrete gradient is defined by:
\be \ba{l} \dis \label{grad-incre} 
 \frac{\Delta H}{\Delta y^1} = \frac{H(y_{n+1}^1, y_n^2, y_n^3,\ldots, y_n^{2 m}) - H (y_n^1, y_n^2, y_n^3,\ldots, y_n^{2 m})}{y_{n+1}^1 - y_n^1}  , \\[3ex]\dis  
 \frac{\Delta H}{\Delta y^2} = \frac{H(y_{n+1}^1, y_{n+1}^2, y_n^3,\ldots, y_n^{2 m}) - H (y_{n+1}^1, y_n^2, \ldots, y_n^{2 m})}{y_{n+1}^2 - y_n^2}  , \\[3ex]\dis  
\dotfill  \  \\[3ex]\dis  
\frac{\Delta H}{\Delta y^{2m}} = \frac{H(y_{n+1}^1, y_{n+1}^2, \ldots, y_{n+1}^{2m}) - H (y_{n+1}^1, y_{n+1}^2, \ldots, y_n^{2m}) }{y_{n+1}^{2m} - y_n^{2m}} .
\ea \ee
where, to fix our attention, $\pmb{y}$ is defined by \rf{y=xp}. The corresponding numerical scheme \rf{disgrad} will be called  GR-IA. In fact, we may identify with $\pmb{y}$ any permutation of $2 m$ components $x^k, p^j$. Thus we have $(2 m)!$ discrete gradients of this type (some of them may happen to be identical). 

Having any discrete gradient we can easily obtain the related symmetric discrete gradient
\be  \label{grad-sym}
 {\bar \nabla}_s H (\pmb{y}_n, \pmb{y}_{n+1}) = \frac{1}{2} \left( {\bar \nabla} H (\pmb{y}_n, \pmb{y}_{n+1}) + {\bar \nabla} H (\pmb{y}_{n+1}, \pmb{y}_n) \right) \ .
\ee
One can easily verify that ${\bar \nabla}_s H$ satisfies conditions \rf{disgrad-def} provided that they are satisfied by $\bar \nabla H$.  Discrete gradient scheme obtained by symmetrization \rf{grad-sym} from GR-IA will be called GR-SYM.

\subsection{Linearization of discrete gradients} 

In order to construct locally exact modifications we need to linearize  discrete gradients defined by \rf{grad-incre} and \rf{grad-sym}. 

\begin{lem} \label{lem-inc}
Linearization of the coordinate increment discrete gradient \rf{grad-incre} around $\pmb{\bar y}$ yields
\be  \label{Tay-AB} 
  \bar \nabla H (\pmb{y}_n, \pmb{y}_{n+1}) \approx H_{\pmb{y}} + \frac{1}{2} \left( A \pmb{\nu}_{n+1} + B \pmb{\nu}_n \right)  , 
\ee
where we denoted $\pmb{\nu}_n = \pmb{y}_n - \pmb{\bar y}$, and 
\be \ba{c} \dis  \label{AB}
A = \left(  \ba{lllll} \frac{1}{2} H_{y^1 y^1} & 0 & \ldots & 0 & 0 \\ H_{y^2 y^1} & \frac{1}{2} H_{y^2 y^2} & \ldots  & 0 & 0 \\ \multicolumn{5}{c}{\dotfill} \\ H_{y^{2 m-1} y^1} & H_{y^{2m-1} y^2} & \ldots & \frac{1}{2} H_{y^{2m-1} y^{2m-1}} & 0 \\ H_{y^{2m} y^1} & H_{y^{2m} y^2} & \ldots & H_{y^{2m} y^{2m-1} } & \frac{1}{2} H_{y^{2m} y^{2m}}   \ea              \right) , 
\\[9ex] \dis
B = \left(  \ba{lllll} \frac{1}{2} H_{y^1 y^1} & H_{y^1 y^2} & \ldots & H_{y^1 y^{2m-1}} & H_{y^1 y^{2m}} \\ 0 & \frac{1}{2} H_{y^2 y^2} & \ldots  & H_{y^2 y^{2m-1}} & H_{y^2 y^{2m}} \\ \multicolumn{5}{c}{\dotfill} \\ 0 & 0 & \ldots & \frac{1}{2} H_{y^{2m-1} y^{2m-1}} & H_{y^{2m-1} y^{2m}} \\ 0 & 0 & \ldots & 0 & \frac{1}{2} H_{y^{2m} y^{2m}}   \ea              \right) , 
\ea \ee
where $\dis H_{y^j y^k} := \frac{\partial^2 H}{\partial y^j   \partial y^k}$. 
\end{lem}

\begin{Proof}
We  denote $\pmb{\hat y}_n^{j}= (y_{n+1}^1,\ldots,y_{n+1}^{j},y_n^{j+1},\ldots, y_n^{2 m})^T$. In particular, $\pmb{\hat y}_n^{0} = \pmb{y}_n$ and $\pmb{\hat y}_n^{2m} = \pmb{y}_{n+1}$. Expanding $H (\pmb{\hat y}_n^{j})$ around $\pmb{\hat y}_n^{j-1}$, we get
\be
 H (\pmb{\hat y}_n^{j}) = H (\pmb{\hat y}_n^{j-1}) + H_{y^j} (\pmb{\hat y}_n^{j-1}) (y_{n+1}^j - y_n^j) + \frac{1}{2} H_{y^j y^j} (\pmb{\hat y}_n^{j-1}) (y_{n+1}^j - y_n^j)^2 + \ldots  
\ee  
Hence
\be  \label{expa1}
\frac{\Delta H}{\Delta y^j} = \frac{ H (\pmb{\hat y}_n^{j}) - H (\pmb{\hat y}_n^{j-1})}{y_{n+1}^j - y_n^j} = 
H_{y^j} (\pmb{\hat y}_n^{j-1}) + \frac{1}{2} H_{y^j y^j} (\pmb{\hat y}_n^{j-1}) (y_{n+1}^j - y_n^j) + \ldots  
\ee
Then, from the definition of $\pmb{\nu}_n$, we have
\be \ba{l} \dis
y_{n+1}^j - y_n^j = \nu_{n+1}^j - \nu_n^j \ , \\[2ex]\dis
\pmb{\hat y}_n^{j-1} = \pmb{\bar y} + \left(  \nu_{n+1}^1,\ldots,\nu_{n+1}^{j-1},\nu_n^j,\ldots,\nu_n^{2 m} \right) \ ,
\ea \ee 
and, taking it into account, we rewrite \rf{expa1} as 
\be
\frac{\Delta H}{\Delta y^j} = H_{y^j} (\pmb{\bar y}) + \sum_{k=1}^{j-1} H_{y^j y^k} (\pmb{\bar y}) \ \nu_{n+1}^k + \sum_{k=j}^{2 m}  H_{y^j y^k} (\pmb{\bar y}) \ \nu_{n}^k + \frac{1}{2} H_{y^j y^j} (\pmb{\bar y}) (\nu_{n+1}^j - \nu_n^j) + \ldots 
\ee
which is equivalent to \rf{Tay-AB}, \rf{AB}. 
\end{Proof}

\begin{lem} \label{lem-sym}
Linearization of the symmetric discrete gradient yields
\be
\bar \nabla_s H (\pmb{y}_n, \pmb{y}_{n+1}) \approx  H_{\pmb{y}} + \frac{1}{2} H_{\pmb{y} \pmb{y}} \left( \pmb{\nu}_n + \pmb{\nu}_{n+1} \right) \ , 
\ee
where $H_{\pmb{y} \pmb{y}} $ is the Hessian matrix of $H$, evaluated at $\pmb{\bar y}$.  
\end{lem}

\begin{Proof}
We observe that  $A + B = H_{\pmb{y} \pmb{y}}$ 
and then we use \rf{grad-sym} and \rf{Tay-AB}.
\end{Proof}

\subsection{Conservative properties of modified discrete gradients}

In the one-dimensional case  locally exact modifications are  clearly energy-preserving, compare \cite{CR-PRE,CR-BIT}. In the general case, conservative properties are less obvious. In this section we present  several useful results.

\begin{lem}  \label{lem-sigma}
We assume that a $2 m \times 2 m$ matrix $\Lambda$ (depending on $h$ and, possibly, on other variables) is skew-symmetric (i.e., $\Lambda^T = - \Lambda$) and 
\be  \label{AS} 
 \lim_{h \rightarrow 0} \frac{\Lambda}{h} = S \ , \qquad S = \mm 0 & I \\ - I & 0 \ema \ . 
\ee
Then, the numerical scheme  
\be  \label{num-A}
  \pmb{y}_{n+1} - \pmb{y}_n = \Lambda {\bar \nabla} H \ ,
\ee
where $\pmb{y}_n$ is defined by \rf{y=xp} and ${\bar \nabla} H$ satisfies \rf{disgrad-def}, is a consistent integrator for \rf{multiham}  preserving the energy integral up to round-off error. 
\end{lem}

\begin{Proof}
The consistency follows immediately form \rf{AS}. The energy preservation can be shown in the standard way. Using the standard scalar product in $\R^{2 m}$, we multiply both sides of \rf{num-A} by ${\bar \nabla} H$
\be
\scal{{\bar \nabla} H}{ \pmb{y}_{n+1} - \pmb{y}_n } = \scal{{\bar \nabla} H}{ \Lambda {\bar \nabla} H } \ .
\ee
By virtue of \rf{disgrad-def} the left-hand side equals $H (\pmb{y}_{n+1}) - H (\pmb{y}_n)$. The right hand side vanishes due to the skew symmetry of $\Lambda$. Hence $H (\pmb{y}_{n+1}) = H (\pmb{y}_n)$. 
\end{Proof}

\begin{rem}
Scheme \rf{num-A} for $\Lambda = h S$  becomes  a standard discrete gradient scheme. 
\end{rem}

\begin{lem}  \label{lem-zachow}
We assume that $2 m \times 2 m$ matrix \ $\pmb{\theta}$ \ is of the following form
\be  \label{theta}
  \pmb{\theta} = \mm \pmb{\delta} &  - \pmb{\sigma} \\ \pmb{\rho} & \pmb{\delta}^T \ema \ , \quad 
\pmb{\rho}^T = - \pmb{\rho} \ , \quad \pmb{\sigma}^T = - \pmb{\sigma} \ , \quad \lim_{h\rightarrow 0} \frac{\pmb{\theta}}{h} = I 
\ee 
(where $\pmb{\delta}$, $\pmb{\sigma}$, $\pmb{\rho}$ are $m \times m$ matrices) and \ ${\bar \nabla} H $ \ is (any) discrete gradient. Then, the numerical scheme given by
\be  \label{mod-locex}
   \pmb{y}_{n+1} - \pmb{y}_n  = \pmb{\theta} S {\bar \nabla} H  
\ee
preserves the energy integral exactly, i.e., $H (\pmb{x}_{n+1}, \pmb{y}_{n+1}) = H (\pmb{x}_n, \pmb{y}_n)$.  
\end{lem}

\begin{Proof} We observe that \ $\pmb{\theta} S$ \ is skew-symmetric, and then we use Lemma~\ref{lem-sigma}.  
\end{Proof}

We easily see that \rf{mod-locex} reduces to the standard discrete gradient scheme when we take \ $\pmb{\theta} = h_n I  $ (i.e., \ $\pmb{\theta}$ is proportional to the unit matrix).

\begin{lem}  \label{lem-anal}
If \ $\pmb{\theta}$ is of the form \rf{theta} and $z \mapsto f (z)$ is any analytic function, then $f (\pmb{\theta})$ is also of the form \rf{theta}.    
\end{lem}

\begin{Proof}
First, we will show that the conditions \rf{theta}  are  equivalent to
\be  \label{theta1}
   \pmb{\theta}^T = S^{-1} \pmb{\theta} S \ . 
\ee
Indeed, assuming  a general form of $\pmb{\theta}$, e.g., $\pmb{\theta} = \mm \pmb{\delta} & - \pmb{\sigma} \\ \pmb{\rho} & \pmb{\gamma} \ema$ we see that  
the constraint \rf{theta1} is equivalent to $\pmb{\gamma} = \pmb{\delta}^T$, $\pmb{\rho}^T = - \pmb{\rho}$ and $\pmb{\sigma}^T = - \pmb{\sigma}$. Then the proof is straightforward. Assuming  $\dis f (z) = \sum_{k=1}^\infty a_n z^n$, we obtain
\be
 \left( \sum_{k=1}^\infty a_k \pmb{\theta}^k \right)^T = \sum_{k=1}^\infty a_k \left( S^{-1} \pmb{\theta} S \right)^k = S^{-1} \left( \sum_{k=1}^\infty a_k \pmb{\theta}^k \right) S \ , 
\ee
i.e., $ f (\pmb{\theta})^T = S^{-1} f (\pmb{\theta}) S$.  
\end{Proof}

\begin{cor}  \label{cor-pres}
If \ $\pmb{\theta}^T = S^{-1} \pmb{\theta} S$ \ and $f$ is an analytic function, then the scheme $\pmb{y}_{n+1} - \pmb{y}_n  = f (\pmb{\theta}) S {\bar \nabla} H $ exactly preserves the energy integral $H$.
\end{cor}

\subsection{Locally exact symmetric discrete gradient scheme}

We begin with the symmetric case because the coordinate increment discrete gradient case is more difficult.  In the symmetric case a locally exact modification is derived similarly as in the one-dimensional case. 

\begin{prop}  \label{prop-sym1-locex}
The following modification of  a  symmetric discrete gradient scheme  is locally exact at \ $\pmb{\bar y}$:  
\be  \label{sym1-locex} 
\pmb{y}_{n+1} - \pmb{y}_n = \pmb{\theta}_n S {\bar \nabla}_s H \ , 
\ee
where 
\be  \label{theta-n}
 \pmb{\theta}_n = 2 (F')^{-1} \tanh \frac{h_n F'}{2} \ , 
\ee
and $F'$, given by \rf{F-F'}, is evaluated at $\pmb{y} = \pmb{\bar y}$. 
\end{prop}

\begin{Proof}  We are going to derive \rf{theta-n}, assuming that $\pmb{\theta}_n$ depends on $\pmb{\bar y}$ and $h$. By virtue of Lemma~\ref{lem-sym} the linearization of \rf{sym1-locex} is given by
\be \label{lin1} 
\pmb{\nu}_{n+1} - \pmb{\nu}_n = \pmb{\theta}_n S \left( H_{\pmb{y}} + \frac{1}{2} H_{\pmb{y} \pmb{y}} (\pmb{\nu}_{n+1} + \pmb{\nu}_n) \right) \ .
\ee
Taking into account \rf{F-F'} we transform \rf{lin1} into 
\be \label{nu-dis} 
 \left( I - \frac{1}{2} \pmb{\theta}_n F' \right) \pmb{\nu}_{n+1} =  \left( I + \frac{1}{2} \pmb{\theta}_n F' \right) \pmb{\nu}_n + \pmb{\theta}_n F \ .
\ee
The scheme \rf{sym1-locex} is locally exact if and only if \rf{nu-dis} coincides with \rf{lin-hex}. Therefore, we require that 
\be   \label{two-eq1}
\left( I - \frac{1}{2} \pmb{\theta}_n F' \right) e^{h_n F'} =  I + \frac{1}{2} \pmb{\theta}_n F'  \ , 
\ee
\be  \label{two-eq2}
\left( I - \frac{1}{2} \pmb{\theta}_n F' \right) \left( e^{h_n F'} - I \right) (F')^{-1} F = \pmb{\theta}_n F \ .
\ee
From equation \rf{two-eq1} we can compute $\pmb{\theta}_n$ which yields \rf{theta-n}. Equation \rf{two-eq2} is automatically satisfied provided that \rf{two-eq1} holds. 
\end{Proof}

Taking the symmetrization of the Itoh-Abe discrete gradient (compare \rf{grad-sym}), modifying it according to \rf{sym1-locex} and choosing $\pmb{\bar y} = \pmb{y}_n$ we get a scheme called GR-SYM-LEX, while for   $\pmb{\bar y} = \frac{1}{2} \left(  \pmb{y}_n + \pmb{y}_{n+1} \right)$  we obtain GR-SYM-SLEX. 

\begin{prop}  \label{prop-sym1-pres}
The numerical scheme \rf{sym1-locex} with $\pmb{\theta}_n$ given by \rf{theta-n} is energy-preserving. 
\end{prop}

\begin{Proof} We have 
$F' = S H_{\pmb{y} \pmb{y}}$. Therefore, $(F')^T = - H_{\pmb{y} \pmb{y}} S = - S^{-1} F' S$, and, as a consequence 
\be  \label{F'2}
((F')^2 )^T = S^{-1} (F')^2 S \ . 
\ee
It means that $(F')^2$ has the form \rf{theta}. The formula \rf{theta-n} expresses $\pmb{\theta}_n$ as an analytic function of $(F')^2$.  Finally, we use Lemma~\ref{lem-anal}. 
\end{Proof}

In the one-dimensional case $(F')^2$ is proportional to the unit matrix which essentially simplifies arguments presented in this section, see \cite{CR-PRE,CR-BIT}.

\subsection{Locally exact coordinate increment discrete gradient scheme}

The symmetric form of the discrete gradient leads to a simple form of the locally exact modification. It turns out, however, that starting from the simplest form of the discrete gradient, namely coordinate increment discrete gradient GR-IA, we also succeed in deriving  corresponding locally exact modifications.

\begin{prop}
The following modification of the coordinate increment discrete gradient scheme is locally exact at \ $\pmb{\bar y}$:  
\be  \label{incre1-locex} 
\pmb{y}_{n+1} - \pmb{y}_n = \pmb{\theta}_n S {\bar \nabla} H \ , 
\ee
where ${\bar \nabla} H$ is given by \rf{grad-incre}, 
\be  \label{theta2-n}
 \theta_n = 2 \left( S R + F' \coth \frac{h_n F'}{2} \right)^{-1} \ , 
\ee
$F'$ is given by \rf{F-F'} (i.e., $F' = S H_{\pmb{y} \pmb{y}}$), and, finally $R = A - B$, i.e., 
\be
 R = \left(  \ba{ccccc} 0  & - H_{y^1 y^2} & \ldots & - H_{y^1 y^{2m -1}} & - H_{y^1 y^{2 m}} \\ H_{y^2 y^1} & 0 & \ldots  & - H_{y^2 y^{2m -1}} & - H_{y^2 y^{2 m}} \\ \multicolumn{5}{c}{\dotfill} \\ H_{y^{2 m-1} y^1} & H_{y^{2m-1} y^2} & \ldots & 0 & - H_{y^{2m-1} y^{2m}}  \\ H_{y^{2m} y^1} & H_{y^{2m} y^2} & \ldots & H_{y^{2m} y^{2m-1} } & 0    \ea    \right)  .     
\ee
$F'$ and $R$ are evaluated at \ $\pmb{y} = \pmb{\bar y}$. 
\end{prop}

\begin{Proof} We are going to derive \rf{theta2-n}, assuming that $\pmb{\theta}_n$ depends on $\pmb{\bar y}$ and $h$. By virtue of Lemma~\ref{lem-inc} the linearization of \rf{incre1-locex} is given by
\be  \label{lin2}
 \pmb{\nu}_{n+1} - \pmb{\nu}_n = \pmb{\theta}_n S ( A \pmb{\nu}_{n+1} + B \pmb{\nu}_n ) + \pmb{\theta}_n  S H_{\pmb{y}} \ .
\ee
Hence, taking into account that $S H_{\pmb{y}} = F$,   
\be \label{nu-dis2} 
 \left( I -  \pmb{\theta}_n  S A \right) \pmb{\nu}_{n+1} =  \left( I +  \pmb{\theta}_n  S B \right) \pmb{\nu}_n + \pmb{\theta}_n F \ .
\ee
The scheme \rf{incre1-locex} is locally exact if and only if \rf{nu-dis2} coincides with \rf{lin-hex}. Therefore, we require that 
\be   \label{2-eq1}
\left( I - \pmb{\theta}_n  S A \right) e^{h_n F'} =  I +  \pmb{\theta}_n  S B  \ , 
\ee
\be  \label{2-eq2}
\left( I -  \pmb{\theta}_n  S A \right) \left( e^{h_n F'} - I \right) (F')^{-1} F = \pmb{\theta}_n F \ .
\ee
Inserting \rf{2-eq1} into \rf{2-eq2} we get
\be
\pmb{\theta}_n  S ( B  + A ) (F')^{-1} F = \pmb{\theta}_n F \ ,
\ee
which is identically satisfied by virtue of 
 $A + B = H_{\pmb{y} \pmb{y}} = S^{-1} F'$,  compare \rf{AB}. The remaining equation, \rf{2-eq1}, defines $\pmb{\theta}_n$: 
\be  \label{defteta}
\pmb{\theta}_n  \left( S A   e^{h_n F'}  + S B \right) =  e^{h_n F'} - I \ .
\ee
In order to simplify \rf{defteta} we introduce $R = A - B$ ($R$ is antisymmetric because $B = A^T$). Then, taking into account $SA+SB= F'$, we get
\be
 S A  = \frac{1}{2}  F'  + \frac{1}{2} S R \ , \qquad S B = \frac{1}{2}  F'  - \frac{1}{2} S R \ . 
\ee
Substituting it into \rf{defteta} we complete the proof. 
\end{Proof}

Choosing  $\pmb{\bar y} = \pmb{y}_n$ we obtain  GR-IA-LEX, while for   $\pmb{\bar y} = \frac{1}{2} \left(  \pmb{y}_n + \pmb{y}_{n+1} \right)$  we get  GR-IA-SLEX.

\begin{prop}
The numerical scheme \rf{incre1-locex} with $\pmb{\theta}_n$ given by \rf{theta2-n} is energy-preserving, i.e., $H (\pmb{y}_{n+1}) = H (\pmb{y}_n)$. 
\end{prop}

\begin{Proof} We have 
\be 
  \pmb{\theta}_n^T = 2 \left( (S R)^T + \left( F' \coth \frac{h_n F'}{2} \right)^T \right)^{-1} = S^{-1}  \pmb{\theta}_n S \ ,
\ee
because 
\be
   (S R)^T = (- R)(-S) = S^{-1} \left( S R \right) S 
\ee
and,  by virtue of Lemma~\ref{lem-anal}, 
\be
\left( F' \coth \frac{h_n F'}{2} \right)^T = S^{-1} \left( F' \coth \frac{h_n F'}{2} \right) S \ , 
\ee 
where we took into account \rf{F'2}. Then, we use Corollary~\ref{cor-pres}. 
\end{Proof}

\ods

Locally exact modification \rf{incre1-locex}  for  Hamiltonian  \ $H = T (\pmb{p}) + V (\pmb{x})$ \ with two degrees of freedom is determined by
\be
\theta_n = h_n \tanhc \frac{h_n F'}{2} \left(  I +\frac{1}{2} h_n  S R \tanhc \frac{h_n F'}{2} \right)^{-1} ,
\ee
where
\be
R = \left(  \ba{cccc}    0 & - V,_{12} & 0 & 0 \\ V,_{12} & 0 & 0 & 0 \\ 0 & 0 & 0 & - T,_{12} \\ 0 & 0 & T,_{12} & 0   \ea  \right)
\ee
In numerical experiments we use $T (\pmb{p}) = \frac{1}{2} \pmb{p}^2$ and  a potential  depending only on $r =|\pmb{x}|$. In this case $T,_{12} = 0$  and \  $V,_{12} = \frac{x^1 x^2}{r} \left( \frac{1}{r} V,_r \right),_r$. After straightforward calculations we reduce \rf{incre1-locex} to:
\be  \ba{l}  \dis
\pmb{x}_{n+1} - \pmb{x}_n  = h D \ \frac{\pmb{p}_{n+1} +\pmb{p}_n}{2}  ,  \\[1ex]  \dis
\pmb{p}_{n+1} - \pmb{p}_n = - h D \ \bar\nabla V - \frac{1}{2} h^2   (\det D) \  V,_{12} \mm 0 & 1 \\ -1 & 0 \ema   ,
\ea \ee
where $\Omega^2 = V,_{\pmb{x}\pmb{x}}$ and $D =  \tanhc \frac{h \Omega}{2}  $.

\section{Order of considered numerical schemes}
\label{sec-prop}

First of all, we recall the Taylor expansion of the exact solution of equation $\pmb{\dot x} = F(\pmb{ x})$ (compare, e.g., \cite{Bu,HLW}): 
\be  \ba{l}  \dis  \label{xth}
\pmb{x} (t+h) = \pmb{x} + h F  + \frac{1}{2} h^2  F' F + \frac{1}{3!} h^3 \left(   F'' (F,F) +(F')^2 F   \right)  + \ldots 
\ea \ee
where all terms on the right-hand side are evaluated at $t$. Note that $F$ is a vector, $F'$ is a matrix, $F''$ is a vector-valued bilinear form, etc.  
Then, $\pmb{x} (t)$ and $\pmb{x} (t + h)$ will be identified with $\pmb{x}_n$ and $\pmb{x}_{n+1}$, respectively. To simplify notation, in this section we replace $h_n$ by $h$.  Computing the order of numerical schemes we use expansions
\be  \label{FFy}  \ba{l} \dis 
F(\pmb{x}_{n+1}) = F + h F' b_1 + h^2 \left(   F' b_2 + \frac{1}{2!} F'' (b_1, b_1)  \right) \\[1em] 
\dis\quad \qquad  + h^3 \left(  F' b_3 +  F'' (b_1,b_2)  +  \frac{1}{3!} F''' (b_1,b_1,b_1) \right) + \ldots 
\\[2em]
 \dis 
F\left(\frac{\pmb{x}_n+\pmb{x}_{n+1}}{2}\right) = F +\frac{1}{2}  h F' b_1 + h^2 \left(   \frac{1}{2} F' b_2 + \frac{1}{8} F'' (b_1, b_1)  \right) \\[1em] 
\dis\quad \qquad  + h^3 \left( \frac{1}{2} F' b_3 + \frac{1}{4} F'' (b_1,b_2)  +  \frac{1}{48} F''' (b_1,b_1,b_1) \right) + \ldots 
\ea \ee
where $b_k$ are coefficients of the Taylor expansion of $\Phi (\pmb{x}_n, h)$, i.e.,  
\be  \label{xb}
\pmb{x}_{n+1} = \pmb{x}_n + b_1 h + b_2 h^2 + b_3 h^3 + \ldots 
\ee
In the case of locally exact schemes we use also the following Taylor expansions
\be  \ba{l}  \dis
\varphi_1 (h F' )= I + \frac{1}{2} h F' + \frac{1}{3!} h^2 (F')^2 + \frac{1}{4!} h^3 (F')^3 + \ldots 
\\[1em]  \dis 
 \tanhc \frac{ h_n F'}{2}  =  I - \frac{1}{12}  h_n^2 (F')^2 + \frac{1}{120}  h_n^4 (F')^4 + \ldots  
\ea \ee
substituting  \rf{FFy} if $F'$ is evaluated at $\pmb{x}_{n+1}$ or $\frac{1}{2} \left( \pmb{x}_n + \pmb{x}_{n+1}   \right)$.  The final result of this analysis (expansions \rf{xb} for particular numerical schemes) is presented below. 

\subsubsection*{First order schemes}
\begin{itemize}
\item   Explicit (forward)  Euler scheme, EEU,  
\be
   \pmb{x}_{n+1} = \pmb{x}_n + h F \  .
\ee
\item   Implicit (backward) Euler scheme,  IEU, 
\be
 \pmb{x}_{n+1} = \pmb{x}_n + h F + h^2 F' F + \ldots  \  
\ee

\end{itemize}

\subsubsection*{Second order schemes}

\begin{itemize}

\item  Locally exact explicit Euler scheme ($\pmb{\bar x} = \pmb{x}_n$), EEU-LEX, 
\be
  \pmb{x}_{n+1} = \pmb{x}_n + h F + \frac{1}{2} h^2 F' F + \frac{1}{6} h^3 (F')^2 F + \ldots 
\ee
\item  Locally exact implicit Euler scheme  ($\pmb{\bar x} = \pmb{x}_n$), IEU-LEX, 
\be
  \pmb{x}_{n+1} = \pmb{x}_n + h F + \frac{1}{2} h^2 F' F  + h^3 \left( \frac{1}{6} (F')^2 F + \frac{1}{2} F'' (F,F) \right) +  \ldots 
\ee
\item   Locally exact implicit Euler scheme  ($\pmb{\bar x} = \pmb{x}_{n+1}$), IEU-ILEX, 
\be
  \pmb{x}_{n+1} = \pmb{x}_n + h F + \frac{1}{2} h^2 F' F  + h^3 \left( \frac{1}{6} (F')^2 F + \frac{1}{4} F'' (F,F) \right) +  \ldots 
\ee
\item   Implicit midpoint rule, IMP, 
\be
  \pmb{x}_{n+1} = \pmb{x}_n + h F + \frac{1}{2} h^2 F' F + h^3 \left( \frac{1}{4} (F')^2 F + \frac{1}{8} F'' (F,F) \right) + \ldots 
\ee
\item   Locally exact implicit midpoint rules IMP-LEX,  IMP-SLEX ($h$-expan\-sions for $\pmb{\bar x} = \pmb{x}_n$ and $\pmb{\bar x} = \frac{1}{2}(\pmb{x}_n + \pmb{x}_{n+1})$ are identical up to the third order):
\be
  \pmb{x}_{n+1} = \pmb{x}_n + h F + \frac{1}{2} h^2 F' F  + h^3 \left( \frac{1}{6} (F')^2 F + \frac{1}{8} F'' (F,F) \right) +  \ldots 
\ee
\item   Trapezoidal rule, TR,
\be
  \pmb{x}_{n+1} = \pmb{x}_n + h F + \frac{1}{2} h^2 F' F + h^3 \left( \frac{1}{4} (F')^2 F + \frac{1}{4} F'' (F,F) \right) + \ldots 
\ee
\item   Locally exact trapezoidal rules TR-LEX, TR-SLEX  ($h$-expansions for $\pmb{\bar x} = \frac{1}{2}(\pmb{x}_n + \pmb{x}_{n+1})$ and  $\pmb{\bar x} = \pmb{x}_n$  are identical up to the third order):
\be
  \pmb{x}_{n+1} = \pmb{x}_n + h F + \frac{1}{2} h^2 F' F  + h^3 \left( \frac{1}{6} (F')^2 F + \frac{1}{4} F'' (F,F) \right) +  \ldots 
\ee

\end{itemize}

\subsection*{Gradient schemes} 

Discussing gradient schemes we confine ourselves to Hamiltonians of the form $H = T (\pmb{p}) + V (\pmb{x}) $. First, we expand a discrete gradient of $V$  with respect to  $\Delta x^j$ \ ($j$th component of $ \pmb{x}_{n+1} - \pmb{x}_n$,  compare \rf{barnablaH}): 
\be  \label{Vexp}
  \frac{\Delta V}{\Delta x^j}  = V,_{x^j}  + \frac{1}{2} A_{jk} \Delta x^k  +  \frac{1}{3!} B_{j\mu\nu} \Delta x^\mu \Delta x^\nu +  \frac{1}{4!} C_{jk\mu\nu} \Delta x^k   \Delta x^\mu   \Delta x^\nu  + \ldots 
\ee
 where  $A_{jk}$, $B_{j\mu\nu}$ and $C_{jk\mu\nu}$ are $\pmb{x}_n$-dependent coefficients. In particular:
\begin{itemize}
\item  Itoh-Abe  gradient ,  GR-IA,  
\be
    A_{jj} = V,_{x^j x^j}  , \quad  A_{jk} = 0 \ \ (j<k) ,  \quad  A_{jk} = 2 V,_{x^j x^k} \ \  (j >k) , 
\ee
\item any symmetric discrete gradient, compare \rf{grad-sym}, 
\be 
 A_{jk} =  V,_{x^j x^k} \ , 
\ee
\item  symmetrization of  the Itoh-Abe gradient, one degree of freedom,
\be
A_{11} =  V,_{xx} \ , \quad B_{111} = V,_{xxx} \ , \quad C_{1111} = V,_{xxxx} \ ,
\ee
\item symmetrization of the Itoh-Abe gradient (GR-SYM), two degrees of freedom,
\be  \ba{l} \dis 
  B_{111} = V,_{xxx} \ , \quad B_{112} = \frac{3}{4} V,_{xxy} \ , \quad B_{122} = \frac{3}{2} V,_{xyy} \ , \\[1em] \dis  B_{222}=V,_{yyy} \ , \quad  B_{212} = \frac{3}{4} V,_{xyy} \ , \quad  B_{211} = \frac{3}{2} V,_{xxy}   \ . 
\ea \ee
\end{itemize}

Having \rf{Vexp}  we substitute expansion of $\Delta x^j$ with respect to $h$ in the considered discrete gradient scheme. As a results we obtain Taylor series for $\pmb{x}_{n+1}$ and comparing them with \rf{xth} we arrive at the following conclusions  (assuming the generic case, because for some exceptional $V$, e.g., $V$ linear or quadratic  in $x^j$, the order can be higher). 

\subsubsection*{First order schemes}

\begin{itemize}

\item  Discrete gradient schemes such that $A_{jk} \neq V,_{x^j x^k}$, in particular: the Itoh-Abe scheme GR-IA.

\end{itemize}

\subsubsection*{Second order schemes}

\begin{itemize}

\item Discrete gradient schemes such that  $A_{jk} = V,_{x^j x^k}$, in particular: symmetric gradient schemes (including GR-SYM) and one-dimensional discrete gradient method.  

\item Locally exact ($\pmb{\bar x}=\pmb{x}_n$, i.e., LEX) modifications of symmetric discrete gradient schemes such that  $B_{j\mu\nu} \neq V,_{x^j x^\mu x^\nu}$. In particular, locally exact modification of symmetrization of the Itoh-Abe discrete gradient scheme (GR-SYM-LEX), multidimensional case (i.e., two degrees of freedom, at least). 

\item Locally exact time reversible ($\pmb{\bar x}=\frac{1}{2} \left( \pmb{x}_n + \pmb{x}_{n+1} \right)$, i.e., SLEX) modifications of symmetric discrete gradient schemes such that  $B_{j\mu\nu} \neq V,_{x^j x^\mu x^\nu}$. 
In particular, locally exact time reversible modification of symmetrization of the Itoh-Abe discrete gradient scheme (GR-SYM-SLEX), multidimensional case.  

\item  Locally exact  ($\pmb{\bar x}=\pmb{x}_n$)  modification of the Itoh-Abe discrete gradient scheme (GR-IA-LEX).

\item  Locally exact time reversible ($\pmb{\bar x}=\frac{1}{2} \left( \pmb{x}_n + \pmb{x}_{n+1} \right)$) modification of  the Itoh-Abe discrete gradient scheme (GR-IA-SLEX).

\end{itemize}

\subsubsection*{Third order schemes}

\begin{itemize}

\item Locally exact (LEX: $\pmb{\bar x}=\pmb{x}_n$) modifications of symmetric discrete gradient schemes such that  $B_{j\mu\nu} = V,_{x^j x^\mu x^\nu}$.  In particular,  GR-LEX in one-dimensional case, see \cite{CR-BIT}. 

\item Locally exact  time reversible ($\pmb{\bar x}=\frac{1}{2} \left( \pmb{x}_n + \pmb{x}_{n+1} \right)$, i.e., SLEX) modifications of symmetric discrete gradient schemes such that   $B_{j\mu\nu} = V,_{x^j x^\mu x^\nu}$ and $C_{j k \mu\nu} \neq V,_{x^j x^k x^\mu x^\nu}$. 

\end{itemize}

\subsubsection*{Fourth  order schemes}

\begin{itemize}

\item Locally exact  time reversible (SLEX: $\pmb{\bar x}=\frac{1}{2} \left( \pmb{x}_n + \pmb{x}_{n+1} \right)$) modifications of symmetric discrete gradient schemes such that   $B_{j\mu\nu} = V,_{x^j x^\mu x^\nu}$ and $C_{j k \mu\nu} =  V,_{x^j x^k x^\mu x^\nu}$.  In particular, GR-SLEX in one-dimensional case, see \cite{CR-BIT}. 

\end{itemize}

\section{Linear stability}
\label{sec-stabil} 

Locally exact numerical schemes have excellent qualitative behaviour around all fixed points of the considered system. 

\begin{prop}
If the equation $\pmb{\dot x} = F (\pmb{x})$ has a fixed point at \ $\pmb{x} = \pmb{\bar x}$, then 
all its locally exact discretizations  have a fixed point at $\pmb{x}_n = \pmb{\bar x}$, as well.  
\end{prop}

\begin{Proof} If \ $F (\pmb{\bar x}) = 0$, then equation \rf{lin-ex} becomes
\be  \label{linex1}
{\pmb x}_{n+1} - \pmb{\bar x} = e^{ h_n F' (\pmb{\bar x})} \left( x_n - \pmb{\bar x} \right) \ . 
\ee
We require that the scheme $\pmb{x}_{n+1} = \Phi (\pmb{x}_n, h)$ is a locally exact discretization of $F (\pmb{\bar x}) = 0$, i.e., its linearization,  given by
\be
\pmb{x}_{n+1} =  \Phi (\pmb{\bar x}, h) + \Phi' (\pmb{\bar x}, h) (\pmb{x}_n - \pmb{\bar x}) \ ,
\ee
coincides with \rf{linex1}. Hence
\be
 e^{ h_n F' (\pmb{\bar x})} = \Phi' (\pmb{\bar x}, h) \ , \qquad \Phi (\pmb{\bar x}, h) = \pmb{\bar x} \ ,
\ee
which means that $\pmb{\bar x}$ is a fixed point of the system  $\pmb{x}_{n+1} = \Phi (\pmb{x}_n, h)$. 
\end{Proof}

Stability of numerical integrators can be roughly defined as follows: 
``the numerical solution provided by a stable numerical integrator does not tend to infinity when the exact solution is bounded'' (see \cite{CCO}, p. 358). The integrator which is stable when applied to linear equations is said to be linearly stable. We may use the notion of A-stability (see, e.g., \cite{Is-book}): an integrator is said to be A-stable, if discretizations of all stable linear equations are stable as well. 

Making one more assumption (quite natural, in fact) that the discretization of a linear system is linear, we see that locally exact integrators are linearly stable. Indeed, solutions of any locally exact discretization have the same trajectories as corresponding exact solutions.  

\begin{cor}
Any locally exact numerical scheme is linearly stable and, in particular, A-stable. 
\end{cor}

What is more, a locally exact discretization yields the best (exact) simulation of a linear equation in the neighbourhood of a fixed point.  In particular, locally exact discretizations preserve any qualitative features of trajectories of linear equations (up to round-off errors, of course).  

Numerical experiments show that locally exact schemes are exceptionally stable also for some simple nonlinear systems \cite{CR-long,CR-PRE,CR-BIT}, but we have no theoretical results concerning the stability (e.g., algebraic stability \cite{Is-book}) in the nonlinear case.

In order to illustrate these general results  we 
will apply four schemes presented above to a linear equation $\pmb{\dot x} = A \pmb{x} $. 
Locally exact explicit Euler scheme and  locally exact implicit Euler scheme yield, respectively,  
\be
  \pmb{x}_{n+1} - \pmb{x}_n = \left( e^{h_n A } - I \right)  \pmb{x}_n \ , 
\ee
\be
  \pmb{x}_{n+1} - \pmb{x}_n = \left( I - e^{- h_n A } - I \right)  \pmb{x}_{n+1} \ . 
\ee
Implicit midpoint and trapezoidal rules yield an identical equation, namely 
\be  \label{ideq} 
 \pmb{x}_{n+1} - \pmb{x}_n = \left( \tanh \frac{h_n A}{2} \right) \left( \pmb{x}_{n+1} + \pmb{x}_n \right) \ .
\ee
Simple calculations show that all resulting equations reduce to the exact discretization: $\pmb{x}_{n+1} = e^{h_n A} \pmb{x}_n$.  
In particular, if  real part of any eigenvalue of $A$ is negative, then $\pmb{x}_n \rightarrow 0$ for $n \rightarrow \infty$.  A-stability is evident.

\section{Numerical experiments}
\label{sec-exp}

In order to illustrate performance of locally exact numerical schemes 
we consider circular orbits for the Hamiltonian
\be
   H (\pmb{x}, \pmb{p}) = \frac{1}{2} |\pmb{p}|^2 + \frac{1}{2} | \pmb{x}|^2 - \frac{1}{30} |\pmb{x}|^3 \ , 
\ee
where the exact solution can be easily found. Namely, initial conditions \ $\pmb{x}_0 = (x_0,0)$, $\pmb{p}_0 = (0, x_0 \sqrt{1 - 0.1 x_0})$ yield a circular orbit  of radius $R=x_0$.
Circular orbits exist for $R<10$. We made numerical tests for  $R=0.2$, $R=1$ and $R=5$ (coresponding periods of exact solutions are given by: $6.347, 6.623$ and $8.886$, respectively).  
While solving implicit algebraic equation for $\pmb{x}_{n+1}$ we used the fixed point method. Iterations were made until the accuracy $3 \cdot 10^{-16}$ was reached. The maximal number of iterations was limited to 20. All figures present  global error at $t=12.5$.

Locally exact modifications improve the accuracy but are quite expensive. In our numerical experiments  we use different time steps for different schemes to assure the same computational cost (i.e., at all figures the computational cost in every column is the same).  The numerical cost is estimated as a number of function evaluations. In other words, we use the scaled time step $\tilde h$ such that $h=\lambda \tilde h$, where $\lambda$  depend on numerical scheme and (to some extent) also on $h$. In order to evaluate the parameter $\lambda$ we start from $\lambda=1$, estimate the computational costs of considered schemes and multiply $h$ by the relative computational cost. The procedure is repeated until computational costs become equal with accuracy about $1 \%$.  Approximate values of the parameter $\lambda$ are given in figure captions.  

\begin{figure} 
\caption{Error of numerical solutions at $t = 12.5$ for a circular orbit ($R=0.2$) as a function of scaled time step $\tilde h$ (i.e.,  $h=\lambda \tilde h$).  EEU: white discs ($\lambda = 1$), IEU: white squares ($\lambda=24,30,42,60$, respectively), EEU-LEX: gray discs ($\lambda=8$), IEU-LEX: gray squares ($\lambda=34,40,53,67$), IEU-ILEX: black squares ($\lambda=130,200,200,200$). }
\label{Euler-0.2}   \par   \ods
\includegraphics[width=\textwidth]{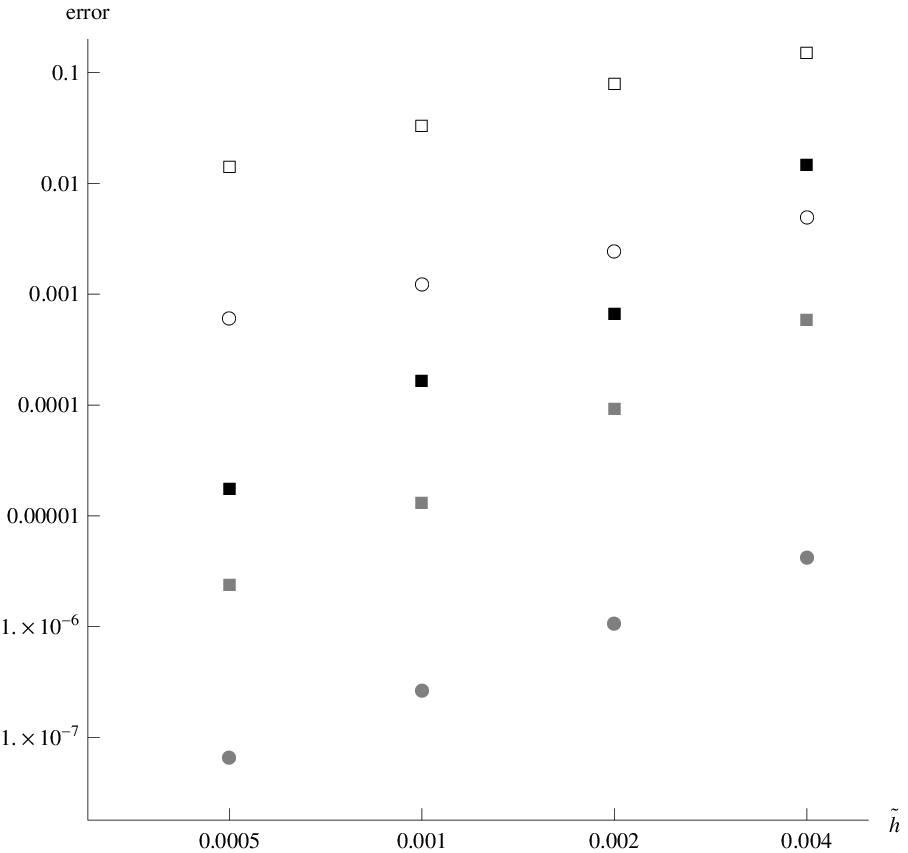} \par
\end{figure}

\begin{figure} 
\caption{Error of numerical solutions at $t = 12.5$ for a circular orbit ($R=5$) as a function of scaled time step $\tilde h$ (i.e.,  $h=\lambda \tilde h$).  EEU: white discs ($\lambda=1$), IEU: white squares ($\lambda=25, 31, 42, 60$), EEU-LEX: gray discs ($\lambda=8$), IEU-LEX: gray squares ($\lambda = 34, 40, 52, 67$), IEU-ILEX: black squares ($\lambda=200$). }
\label{Euler-5}   \par   \ods
\includegraphics[width=\textwidth]{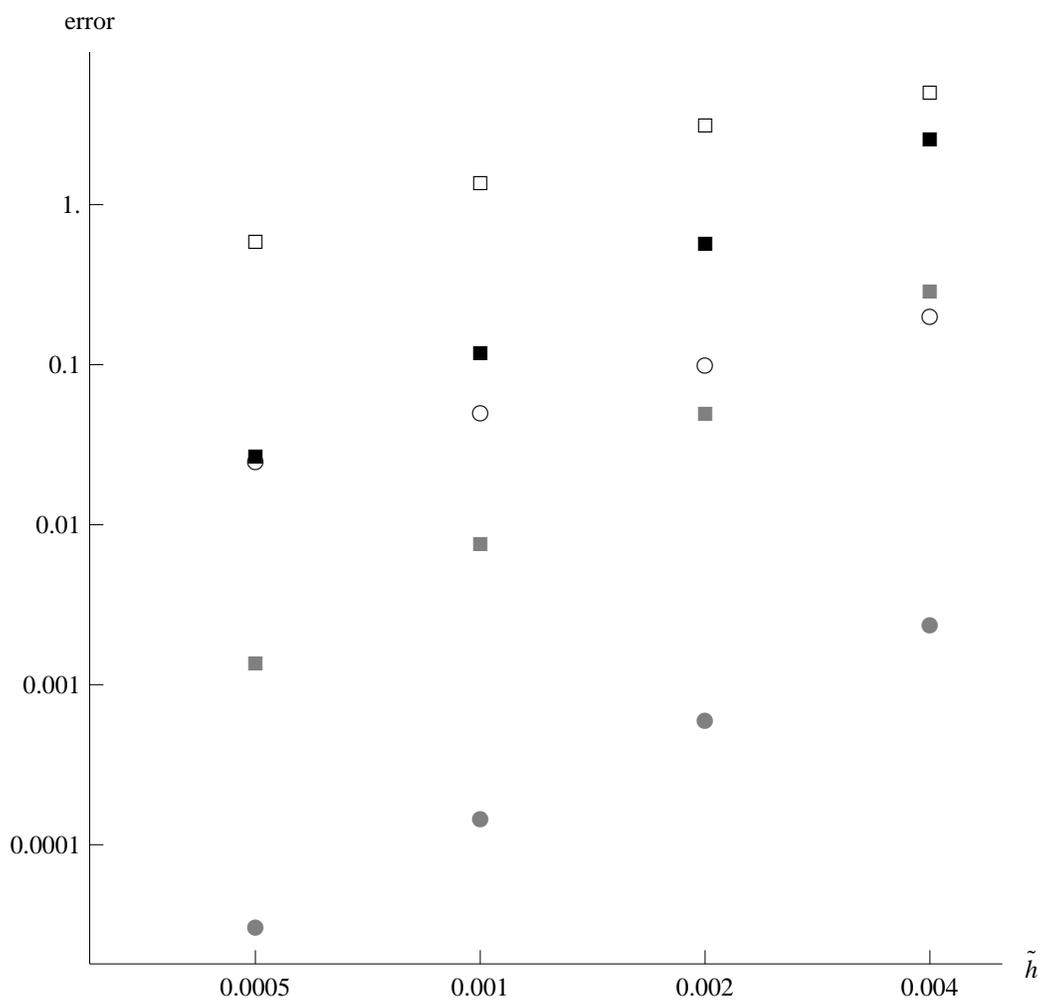} \par
\end{figure}

In the case of Euler schemes  (Figs.~\ref{Euler-0.2} and \ref{Euler-5}) the advantage of locally exact modifications is obvious. Although the cost of locally exact modifications is considerably higher, they perform much better than both standard Euler schemes. Actually the accuracy of all 3 modifications (taken with the same time step) is similar. However, taking into account the computational cost,  we see that the scheme EEU-LEX is the best modification. 

\begin{figure} 
\caption{Error of numerical solutions at $t = 12.5$ for a circular orbit ($R=0.2$) as a function of scaled time step $\tilde h$ (i.e., $h=\lambda \tilde h$).  IMP: white discs ($\lambda=1$), IMP-LEX: gray discs ($\lambda=1.3, 1.3, 1.3, 1.2$), IMP-SLEX: black discs ($\lambda= 3.4, 3.4, 3.9, 3.9$), TR: white squares ($\lambda=1.06$), TR-LEX: gray squares ($\lambda=1.9, 1.6, 1.8, 1.8$), TR-SLEX: black squares ($\lambda=3.5, 3.4, 4.0, 3.9$).}
\label{midtrap-0.2}   \par   \ods
\includegraphics[width=\textwidth]{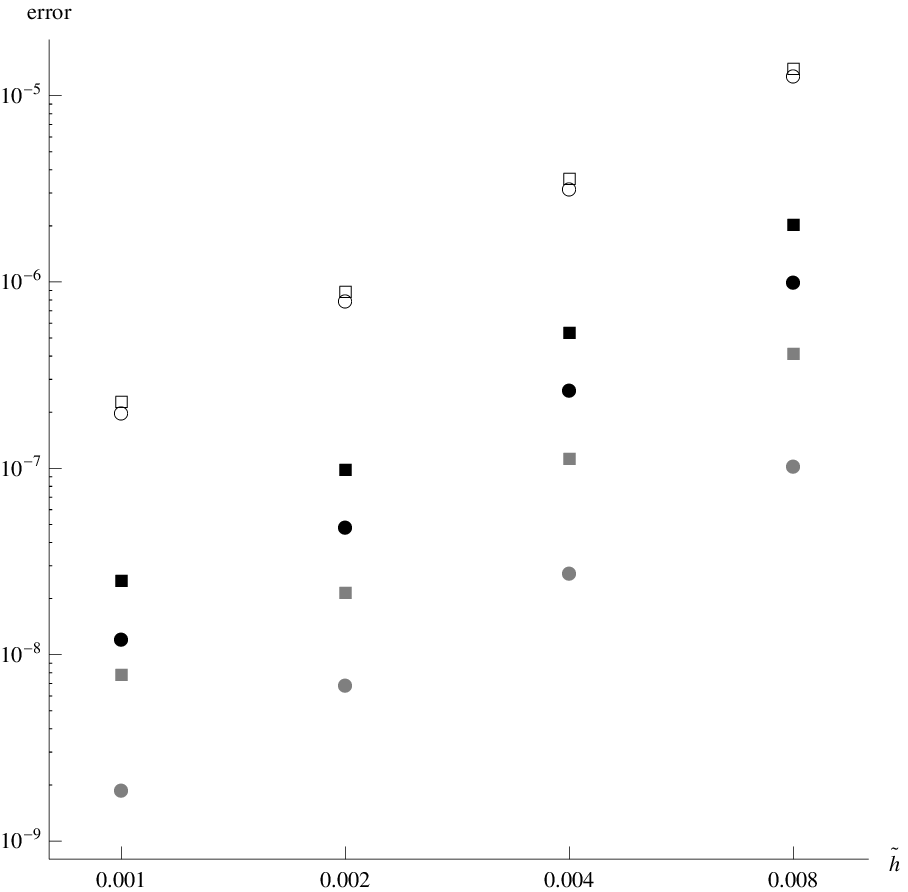} \par
\end{figure}

\begin{figure} 
\caption{Error of numerical solutions at $t = 12.5$ for a circular orbit ($R=1$) as a function of scaled time step $\tilde h$ (i.e., $h=\lambda \tilde h$).  IMP: white discs ($\lambda=1$), IMP-LEX: gray discs ($\lambda=1.5, 1.3, 1.4, 1.2$), IMP-SLEX: black discs ($\lambda= 3.4, 3.4, 3.9, 3.9$), TR: white squares ($\lambda=1.06$), TR-LEX: gray squares ($\lambda=1.9, 1.6, 1.8, 1.8$), TR-SLEX: black squares ($\lambda=3.5, 3.4, 4.0, 3.9$).}
\label{midtrap-1}   \par   \ods
\includegraphics[width=\textwidth]{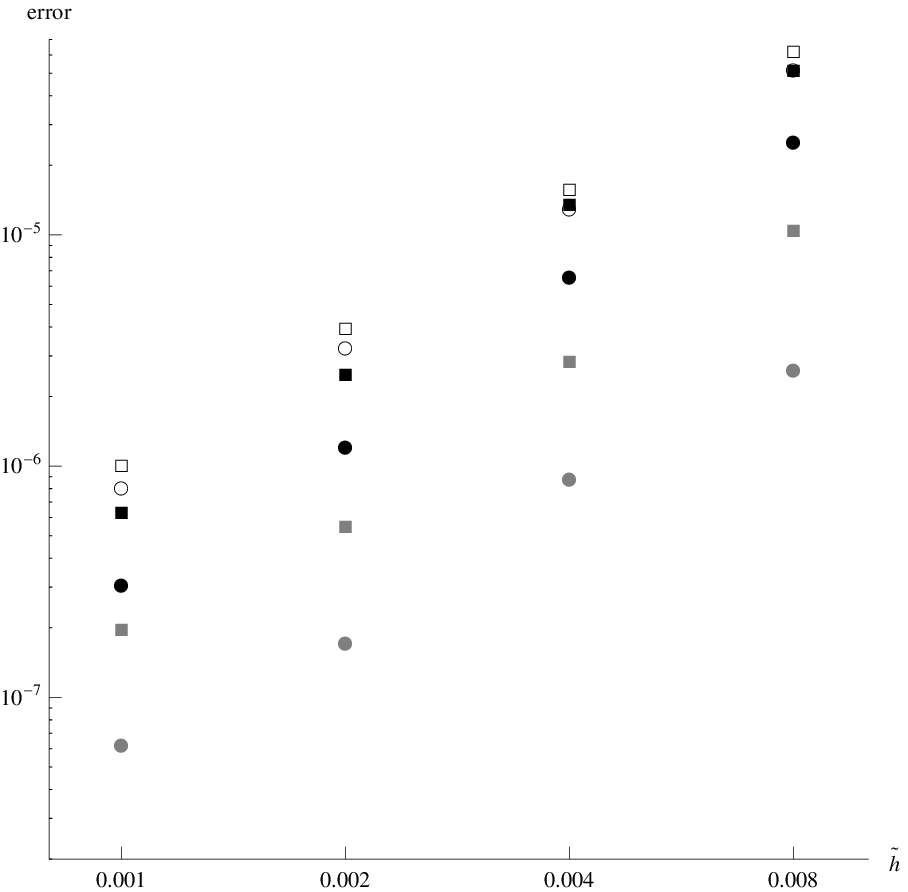} \par
\end{figure}

\begin{figure} 
\caption{Error of numerical solutions at $t = 12.5$ for a circular orbit ($R=0.2$)  as
 a function of scaled time step   $\tilde h$ (i.e., $h=\lambda \tilde h$).  GR-IA: white discs, GR-IA-LEX: gray discs,  GR-IA-SLEX: black discs, GR-SYM: white squares, GR-SYM-LEX: gray squares,   GR-SYM-SLEX: black squares. } 
 \label{grad-0.2}  \par   \ods
 \includegraphics[width=\textwidth]{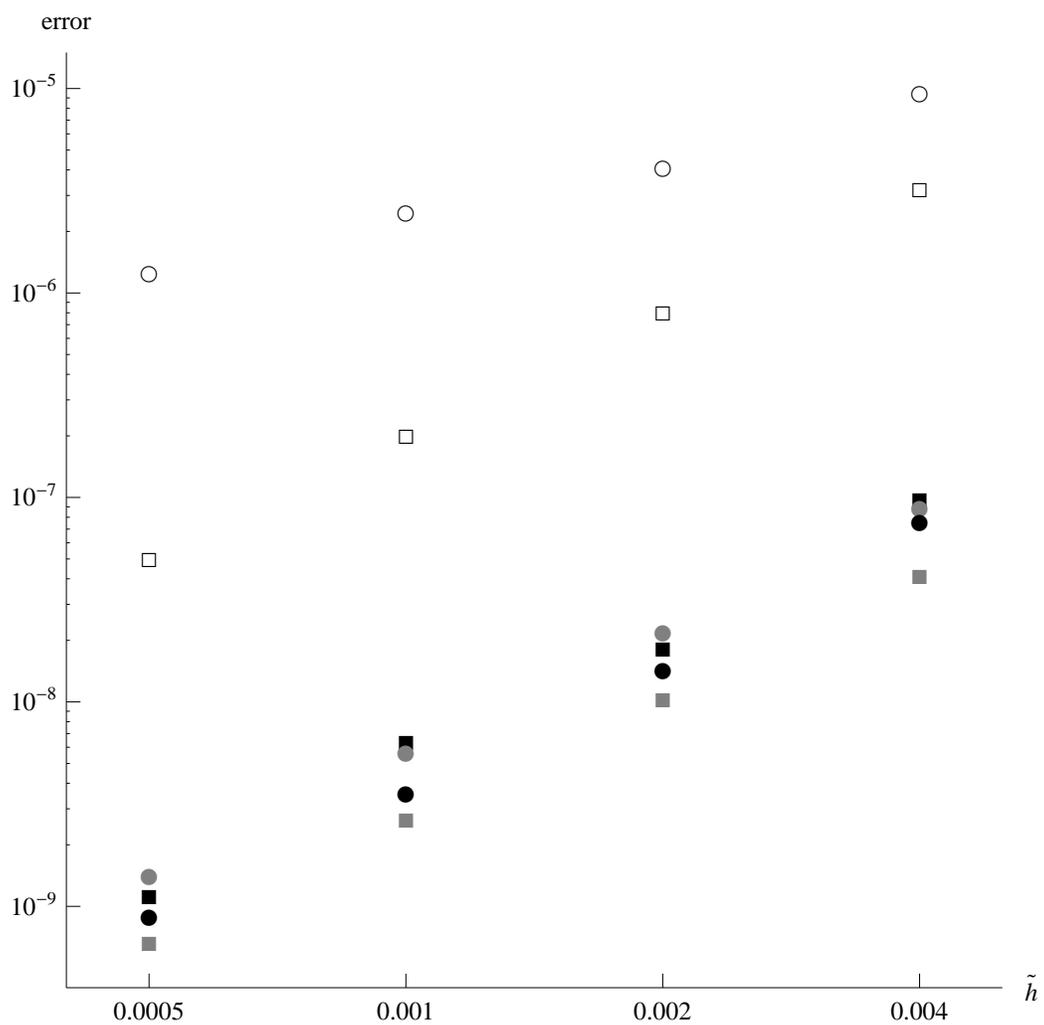} \par
\end{figure}

\begin{figure} 
\caption{Error of numerical solutions at $t = 12.5$ for a circular orbit ($R=1$)  as
 a function of scaled time step $\tilde h$ (i.e., $h=\lambda \tilde h$).  GR-IA: white discs ($\lambda=0.8$), GR-IA-LEX: gray discs ($\lambda=1.1$),  GR-IA-SLEX: black discs ($\lambda=1.5, 1.8, 1.5, 1.7$), GR-SYM: white squares ($\lambda=1$), GR-SYM-LEX: gray squares ($\lambda=1.3, 1.4, 1.3, 1.4$),   GR-SYM-SLEX: black squares ($\lambda=1.7, 2.0, 1.7, 1.9$). } 
 \label{grad-1}  \par   \ods
 \includegraphics[width=\textwidth]{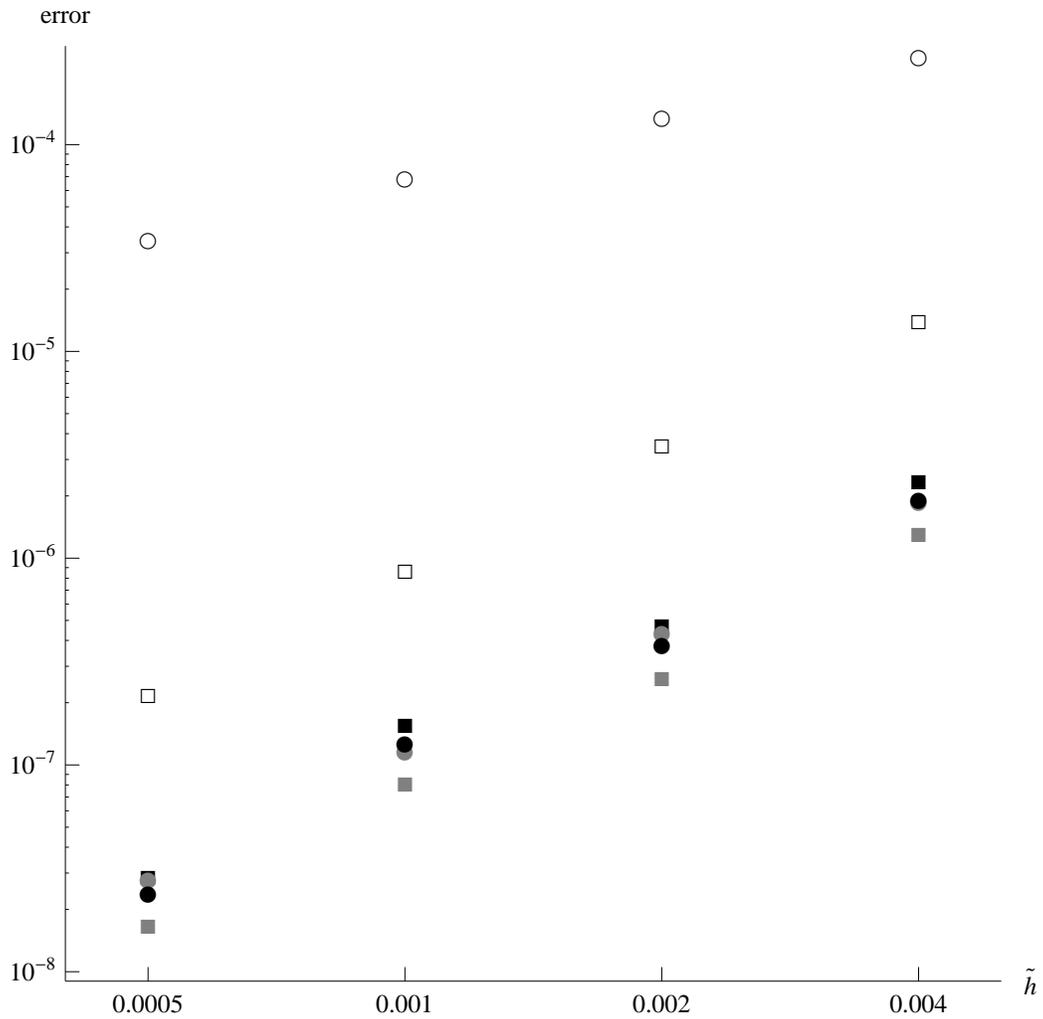} \par
\end{figure}

Locally exact modifications of implicit midpoint and trapezoidal rules increase have higher accuracy only for orbits with a relatively small radius (Figs.~\ref{midtrap-0.2} and \ref{midtrap-1}). The best results are produced by scheme IMP-LEX. For larger orbits the unmodified implicit midpoint rule is most accurate.  Similar situation has place for gradient schemes. Locally exact modifications give a considerable improvement, by about 1-2 orders of magnitude, for 
trajectories in a quite large  neighborhood of the stable equilibrium, see Figs.~\ref{grad-0.2} and \ref{grad-1}. For  large $R$ (e.g., $R=5$)  locally exact modifications improve considerably only GR-IA. Their influence on the GR-SYM is neutral or even negative. In fact GR-SYM has similar accuracy for all considered orbits while its locally exact modifications are very accurate only for smaller values of  $R$.

\section{Concluding remarks}

We presented a new class of numerical schemes characterized by the so called local exactness. This notion is known (under different names) since almost fifty years. The original application, see \cite{Pope}, has been confined to the exact discretization of linearized equations, compare Section~\ref{sec-ex-lin}. We obtain in this way a particular locally exact integrator which has some advantages (e.g., it can serve as a good predictor,  see \cite{CR-PRE}, Section V.A).   

Our approach has two new features. First, we modify known numerical schemes in a locally exact way. Any numerical scheme admits at least one (usually more) natural locally exact modification, see Section~\ref{sec-pop}. Second, we try to preserve geometric properties of the original numerical scheme. This task is not trivial. In this paper we present one successful application: locally exact modifications of discrete gradient methods for canonical Hamilton equations. 
We are able to exactly preserve the energy integral and, in the same time, increase the accuracy by many orders of magnitude. Another advantage is a variable time step. Unlike symplectic methods (which work mostly for the constant time step) discrete gradient methods admit conservative modifications with variable time step. Therefore, one may easily  implement any variable step method in order to obtain further improvement. 

In the one-dimensional case  proposed modifications, although a little bit more expensive, turns out to be more accurate even by 8 orders of magnitude in comparison to the standard discrete gradient scheme, see  \cite{CR-PRE,CR-BIT}.  In multidimensional cases the relative cost of our algorithm is higher, but still  our method  is  of great advantage for orbits in the neighborhood of the stable equilibrium. The accuracy is increased by 1-2 orders of magnitude. We point out that in most cases  locally exact modifications do not change orders of  modified schemes.  

We presented locally exact modifications from a unified theoretical perspective. There are many possible further developments. First of all, we plan to apply our approach to chosen multidimensional problems, testing the accuracy of locally exact schemes by numerical experiments. Then, we would like to extend the range of applications. It would be natural to consider locally exact modifications of Runge-Kutta schemes and $G$-symplectic integrators \cite{Bu}. However, our first attempts seem to suggest that this is a challanging problem. 
Modifications proposed in this paper contain exponentials of variable matrices. Similar time-consuming evaluations are characteristic for all exponential integrators and  in this context 
effective methods of computing matrix exponentials have been recently developed \cite{HoL1,NW}.

Throughout this paper we assumed the autonomous case, ${\pmb {\dot x}} = F ({\pmb x})$. The extension on the non-autonomous case can be done along lines indicated already in Pope's paper \cite{Pope}. A separate problem is to obtain in this case any  locally exact modification with geometric properties.
Another open problem is the construction of locally exact (or, at least, linearization preserving) defomations of generalized discrete gradient algoritms preserving all first integrals (see \cite{MQR2}).  
Linearization-preserving integrators form an important subclass of locally exact schemes. It would be worthwhile to study linearization-preserving modifications of geometric numerical integrators, especially in those cases when locally exact are difficult or impossible to construct.

\ods

{\bf Acknowledgment}. Research supported in part  by the National Science Centre (NCN) grant no. 2011/01/B/ST1/05137.

\end{document}